\documentclass[12pt]{article}
\usepackage{cite}
\usepackage{amssymb}
\usepackage{amsmath}
\usepackage{graphicx}
\newtheorem{thm}{Theorem}[section]
\newtheorem{cor}[thm]{Corollary}

\usepackage{graphics}
\numberwithin{equation}{section}

\newcommand{\eh}{\hfill}\newlength{\sperr}

\newenvironment{proof}{{\settowidth{\sperr}{\bf\rm
Proof}%
\par\addvspace{0.3cm}\noindent\parbox[t]{1.3\sperr}
{\textit{ P\eh r\eh o\eh o\eh f\eh .}}%
}}{\nopagebreak\mbox{}
$\Box$\par\addvspace{0.3cm}}

\def\p{\partial}
\def\ve{\varepsilon}

\newtheorem{Pa}{Paper}[section]
\newtheorem{Tm}[Pa]{{\bf Theorem}}
\newtheorem{La}[Pa]{{\bf Lemma}}
\newtheorem{Cy}[Pa]{{\bf Corollary}}
\newtheorem{Rk}[Pa]{{\bf Remark}}

\newtheorem{Ee}[Pa]{{\bf Example}}
\newtheorem{Dn}[Pa]{{\bf Definition}}
\newtheorem{Pn}[Pa]{{\bf Proposition}}
\numberwithin{equation}{section}

\title{Levy processes: long time behavior and 
convolution-type form of the 
Ito representation
of the  infinitesimal generator}
\author{Lev Sakhnovich}
\date{}
\begin{document}
\maketitle

\thanks{99 Cove ave., Milford, CT, 06461, USA \\
 E-mail: lsakhnovich@gmail.com}\\

 potential,
 quasi-potential, long time behavior.
\begin{abstract}In the present paper we show that the  Levy-Ito representation
of the infinitesimal generator $L$ for Levy processes $X_t$ can be written in a
convolution-type form. Using the obtained convolution form we have constructed
the quasi-potential operator $B$. We denote by $p(t,\Delta)$ the probability
that a sample of the process $X_t$ remains inside the domain $\Delta$ for
 $0{\leq}\tau{\leq}t$ (ruin problem). With the help of the operator $B$ we find
 a new formula for $p(t,\Delta)$. This formula allows us to obtain long time behavior of $p(t,\Delta)$.

  \end{abstract}

\section{Introduction}\label{sec1}
%%%%%%%%

 Let us introduce  the notion of the Levy processes.

\begin{Dn}\label{Definition 1.1.}
A stochastic process $\{X_{t}:t{\geq}0\}$ is called Levy
process, if the following conditions are fulfilled:\\
1. Almost surely $X_{0}=0$, i.e. $P(X_{0}=0)=1$.\\
(One says that an event happens almost surely (a.s.) if it happens with probability one.)\\
2. For any $0{\leq}t_{1}<t_{2}...<t_{n}<\infty$ the random variables \\
 $X_{t_{2}}-X_{t_{1}}, X_{t_{3}}-X_{t_{4}},..., X_{t_{n}}-X_{t_{n-1}}$\\ are independent (independent increments).\\
( To call the increments of the process $X_{t}$ \emph{independent} means that
increments $X_{t_{2}}-X_{t_{1}}, X_{t_{3}}-X_{t_{4}},..., X_{t_{n}}-X_{t_{n-1}}$ are mutually  (not just
pairwise) independent.)\\
3. For any $s<t$ the distributions of   $X_{t}-X_{s}$ and $X_{t-s}$ are equal
(stationary increments).\\
4. Process $X_{t}$ is almost surely right continuous with left limits.\\
\end{Dn}
Then Levy-Khinchine formula gives (see \cite{Bert}, \cite{Sato})
 \begin{equation}
\mu(z,t)=E\{\mathrm{exp}[izX_{t}]\}=
\mathrm{exp}[-t\lambda(z)],\quad t{\geq}0,  \label{1.1} \end{equation}
where
\begin{equation}
\lambda(z)=\frac{1}{2}Az^{2}-i{\gamma}z-\int_{-\infty}^{\infty}(e^{ixz}-1-ixz1_{|x|<1})\nu(dx).
 \label{1.2} \end{equation} Here $A{\geq}0,\quad \gamma=\overline{\gamma},\quad z=\overline{z}$ and
 $\nu(dx)$ is a measure on the axis $(-\infty,\infty)$
  satisfying the conditions
\begin{equation}
\int_{-\infty}^{\infty}\frac{x^{2}}{1+x^{2}}\nu(dx)<\infty.
 \label{1.3} \end{equation}
The  Levy-Khinchine formula is determined by the Levy-Khinchine triplet
(A,$\gamma, \nu(dx)$).\\
By $P_{t}(x_{0},\Delta)$ we denote the probability
$P(X_{t}{\in}\Delta)$ when $P(X_{0}=x_{0})=1$ and $\Delta{\in}R$.
The transition operator $P_{t}$ is defined by the formula
\begin{equation}
P_{t}f(x)=\int_{-\infty}^{\infty}P_{t}(x,dy)f(y). \label{1.4} \end{equation} Let
$C_{0}$ be the Banach space of continuous functions $f(x)$ ,
satisfying the condition $\mathrm{lim}f(x)=0,\quad |x|{\to}\infty$
with the norm $||f||=\mathrm{sup}_{x}|f(x)|$. We denote by
$C_{0}^{n}$ the set of $f(x){\in}C_{0}$ such that
$f^{(k)}(x){\in}C_{0},\quad (1{\leq}k{\leq}n).$ It is known that
\cite{Sato}
\begin{equation} P_{t}f{\in}C_{0}, \label{1.5} \end{equation}
if $f(x){\in}C_{0}^{2}.$\\
Now we formulate the following important result (see \cite{Sato}) .
\begin{Tm}\label{Theorem 1.2.} \textbf{(Levy-Ito decomposition.)}
{The family of the operators $P_{t}\quad
(t{\geq}0)$ defined by the Levy process $X_{t}$ is a strongly
continuous  semigroup on $C_{0}$ with the norm $||P_{t}||=1$. Let
$L$ be its infinitesimal generator. Then}
\begin{equation} Lf=\frac{1}{2}A\frac {d^{2}f}{dx^{2}}+\gamma
\frac{df}{dx}+\int_{-\infty}^{\infty}(f(x+y)-f(x)-y\frac{df}{dx}1_{|y|<1})\nu(dy), \label{1.6} \end{equation}
{where} $f{\in}C_{0}^{2}$.
\end{Tm}
Slightly changing the Sato classification \cite {Sato} we introduce the following definition:\\
\begin{Dn}\label{Definition 1.3} We say that a Levy process $X_{t}$ generated by $(A,\nu,\gamma)$ has  type
$I$ if
\begin{equation}A=0 \,and\,\int_{-\infty}^{\infty}\nu(dx)<\infty,\label{1.7}\end{equation}
and $X_{t}$ has  type $II$ if
\begin{equation}A{\ne}0 \quad or\quad \int_{-\infty}^{\infty}\nu(dx)=\infty.\label{1.8}\end{equation}
\end{Dn}
\begin{Rk}\label{Remark 1.4}The introduced type $I$ coincides with the type $A$ in the
Sato classification.The introduced type $II$ coincides with the union  of the types $B$ and $C$ in the Sato classification.\end{Rk}
The properties of these two types of the Levy processes are quiet different.
2. In the present paper we show that the Levy-Ito representation of the generator $L$ can be written in the convolution form
\begin{equation}
Lf=\frac{d}{dx}S\frac{d}{dx}f, \label{1.9}
\end{equation}
where the operator $S$ is defined by the relation
\begin{equation}
Sf=\frac{1}{2}Af+\int_{-\infty}^{\infty}k(y-x)f(y)dy. \label{1.10} \end{equation}
We note that  for arbitrary $a\, (0<a<\infty)$ the inequality
\begin{equation} \int_{-a}^{a}|k(t)|dt<\infty \label{1.11} \end{equation}
is true.\\ Formulas \eqref{1.9} and \eqref{1.10} were proved  before in our works \cite{Sakh6} under some additional
conditions. In the present paper we omit these additional conditions and prove these formulas for the general case.The representation of $L$ in form \eqref{1.9} is convenient as
the operator $L$ is expressed with the help of the classic
differential and convolution operators. Assuming that $X_t$ belongs to the class II we have constructed the
quasi-potential operator $B.$ This operator $B$ is linear and bounded in the space
of the continuous functions. We denote by $p(t,\Delta)$ the probability  that a sample of the process $X_t$ remains inside
the domain $\Delta$ for $0{\leq}\tau{\leq}t$ (ruin problem). With the help of the operator $B$ we find a new
formula for $p(t,\Delta)$.
This formula allow us to obtain the long time behavior of $p(t,\Delta)$. Namely, we have proved the following asymptotic formula
\begin{equation}
 p(t,\Delta)=e^{-t/\lambda_{1}}[c_{1}+o(1)],\,c_1>0,\,\lambda_{1}>0,\,\quad t{\to}+\infty.
  \label{1.12}
  \end{equation}
Let $T_{\Delta}$  be thus time  during which $X_{t}$ remains in the domain $\Delta$
before it leaves the domain $\Delta$ for the first time. It is easy to see that
\begin{equation}p(t,\Delta)=P(T_{\Delta}>t).\label{1.13}\end{equation}

In Sections \ref{sec1} and \ref{sec2} we often follow the presentation from
\cite{SakhArx}.
\begin{Rk}\label{Remark 1.5}In the next paper we plan to apply the convolution
representation \eqref{1.9} to the Levy processes $X_t$ of the type I.\end{Rk}

%%%%%%
\section{Convolution-type form of infinitesimal generator}
\label{sec2}\paragraph{1.}

By $C(a)$ we denote the
set of functions $f(x){\in}C_{0}$ which have the following property:
 \begin{equation} f(x)=0 ,\quad x{\notin}[-a,a]\label{2.1} \end{equation}
i.e. the function $f(x)$ is equal to zero in the neighborhood of $x=\infty$.
 We note, that parameter $a$
can be different for different $f$.\\
We introduce the functions
\begin{equation}\mu_{-}(x)=\int_{-\infty}^{x}\nu(dx),\,x<0,\label{2.2} \end{equation}
\begin{equation}\mu_{+}(x)=-\int_{x}^{\infty}\nu(dx),\,x>0,\label{2.3} \end{equation}
where the functions $\mu_{-}(x)$ and $\mu_{+}(x)$ are  monotonically increasing and right continuous on the  half-axis   $(-\infty,0]$ and $[0,\infty)$ respectively. We note that
\begin{equation}\mu_{+}(x){\to}0,\,x{\to}+\infty;\,
\mu_{-}(x){\to}0,\,x{\to}-\infty,\label{2.4}\end{equation}
\begin{equation}\mu_{-}(x){\geq}0,\,x<0;\,\mu_{+}(x){\leq}0,\,x>0.\label{2.5}\end{equation} In view of \eqref{1.3}
the integrals in the right sides of \eqref{2.2} and \eqref{2.3} are convergent.
\begin{Tm}\label{Theorem 2.1} The following relations
\begin{equation}{\ve}^{2}\mu_{\pm}(\pm{\ve}){\to}0,\,\ve{\to}+0,\label{2.6} \end{equation}
\begin{equation}\int_{-a}^{0}x\mu_{-}(x)dx<\infty,\,
-\int_{0}^{a}x\mu_{+}(x)dx<\infty,\,0<a<\infty\label{2.7} \end{equation}
are true\end{Tm}
%\emph{Proof.}
\begin{proof} According to \eqref{1.3} we have
\begin{equation}0{\leq}\int_{-a}^{-\ve}x^{2}d\mu_{-}(x){\leq}M,
\label{2.8}\end{equation}where $M$ does not depend from $\ve$.
Integrating by parts the integral of \eqref{2.8} we  obtain:
\begin{equation}\int_{-a}^{-\ve}x^{2}d\mu_{-}(x)=\ve^{2}\mu_{-}(-\ve)-a^{2}\mu_{-}(-a)
-2\int_{-a}^{-\ve}x\mu_{-}(x)dx{\leq}M,\label{2.9}\end{equation} The function $-\int_{-a}^{-\ve}x\mu_{-}(x)dx$ of $\ve$
is monotonic increasing. In view of
\eqref{2.9} this function is bounded. Hence we have
\begin{equation}\lim_{\ve{\to}+0}\int_{-a}^{-\ve}x\mu_{-}(x)dx=\int_{-a}^{0}x\mu_{-}(x)dx
\label{2.10}\end{equation} It follows from \eqref{2.9} and \eqref{2.10} that
\begin{equation}\lim_{\ve{\to}+0}\ve^{2}\mu_{-}(-\ve)=m,\quad m{\geq}0.\label{2.11}\end{equation}
Using  \eqref{2.10} and \eqref{2.11} we have $m=0$.
Thus, relations \eqref{2.9} and \eqref{2.10} are proved for $\mu_{-}(x)$. In the same way
relations \eqref{2.9} and \eqref{2.10} can be proved  for $\mu_{+}(x)$.
\end{proof}
\paragraph{2.}Let us introduce the functions
\begin{equation}k_{-}(x)=\int_{-1}^{x}\mu_{-}(t)dt,\quad -\infty{\leq}x<0,\label{2.11'}
\end{equation}
\begin{equation}k_{+}(x)=-\int_{x}^{1}\mu_{+}(t)dt,\quad 0<x{\leq}+\infty.\label{2.12}
\end{equation}
 In view of \eqref{2.5}
the integrals on the right-hand sides of \eqref{2.11'} and \eqref{2.12} are absolutely convergent.
From \eqref{2.11'} and \eqref{2.12} we obtain the assertions:
\begin{Tm}\label{Theorem 2.2}
1.The function $k_{-}(x)$  is continuous, monotonically increasing on the   $(-\infty,0)$ and
\begin{equation}k_{-}(x){\geq}0,\quad -1{\leq}x<0.\label{2.13}\end{equation}
2.The function $k_{+}(x)$  is continuous, monotonically decreasing on the   $(0,+\infty)$ and
\begin{equation}k_{+}(x){\geq}0,\quad 0<x{\leq}1.\label{2.14}\end{equation}
\end{Tm}
Further we need the following result:
\begin{Tm}\label{Theorem 2.3} The  relations
\begin{equation}{\ve}k_{-}(-\ve){\to}0,\,\, \ve{\to}+0;\quad {\ve}k_{+}(\ve){\to}0,\,\, \ve{\to}+0;
\label{2.15}\end{equation}
\begin{equation}\int_{-1}^{0}k_{-}(x)dx<\infty;\quad \int_{0}^{1}k_{+}(x)dx<\infty
\label{2.16}\end{equation}are valid.
\end{Tm}
%\emph{Proof.}
\begin{proof} According to \eqref{2.5} we have
\begin{equation}0{\leq}-\int_{-1}^{-\ve}x\mu_{-}(x)dx{\leq}M,,
\label{2.18}\end{equation}where $M$ does not depend from $\ve$.
Integrating by parts the integral  of \eqref{2.18} we  obtain:
\begin{equation}-\int_{-1}^{-\ve}x\mu_{-}(x)dx={\ve}k_{-}(-\ve)-k_{-}(-1)
+\int_{-1}^{-\ve}k_{-}(x)dx{\leq}M,\label{2.19}\end{equation}
The function $\int_{-1}^{-\ve}k_{-}(x)dx$ of $\ve$
is monotonic increasing.This function is bounded (see \eqref{2.19}) . Hence we have
\begin{equation}\lim_{\ve{\to}+0}\int_{-1}^{-\ve}k_{-}(x)dx=\int_{-1}^{0}k_{-}(x)dx
\label{2.20}\end{equation} It follows from \eqref{2.19} and \eqref{2.20} that
\begin{equation}\lim_{\ve{\to}+0}{\ve}k_{-}(-\ve)=p,\quad p{\geq}0.\label{2.21}\end{equation}
Using  \eqref{2.20} and \eqref{2.21} we have $p=0$.
Thus, the first relations (i.e., the relations for $k_{-}(x)$) in  \eqref{2.15} and \eqref{2.16}  are proved. 
The second relations (i.e., the relations for $k_{+}(x)$) in \eqref{2.15} and \eqref{2.16}  can be proved  in the same way.
\end{proof}

\paragraph{3.}
We use  the following notation
\begin{equation}J(f)=J_{1}(f)+J_{2}(f),\label{2.22}\end{equation}where
\begin{equation}
J_{1}(f)=\frac{d}{dx}\int_{-\infty}^{x}f^{\prime}(y)k_{-}(y-x)dy,\quad  f(x){\in}C(a), \label{2.23} \end{equation}
\begin{equation}
J_{2}(f)=\frac{d}{dx}\int_{x}^{\infty}f^{\prime}(y)k_{+}(y-x)dy, \quad  f(x){\in}C(a). \label{2.24} \end{equation}
\begin{La}\label{Lemma 2.4.}
The operator $J(f)$ defined by \eqref{2.22} can be represented in the form
\begin{equation}
J(f)=\int_{-\infty}^{\infty}[f(y+x)-f(x)-y\frac{df(x)}{dx}1_{|y|{\leq}1}]{\mu}(dy)
+{\Gamma}f^{\prime}(x),
 \label{2.25} \end{equation}
 where $\Gamma=\overline{\Gamma}$ and
$f(x){\in}C(a)$.
\end{La}

%\emph{Proof.}
\begin{proof}
From \eqref{2.24}  we obtain  the relation
\begin{equation}
J_{1}(f)=-\int_{x-1}^{x}[f^{\prime}(y)-f^{\prime}(x)]
k_{-}^{\prime}(y-x)dy-
\int_{-a}^{x-1}f^{\prime}(y)k_{-}^{\prime}(y-x)dy. \label{2.26} \end{equation}
 By proving \eqref{2.26} we used relations
 \eqref{2.12} and equality
 \begin{equation}\int_{x-1}^{x}k_{-}(y-x)dy=\int_{-1}^{0}k_{-}(v)dv.
 \label{2.27}\end{equation}We introduce the notations \begin{equation}
P_{1}(x,y)=f(y)-f(x)-(y-x)f^{\prime}(x),\quad
P_{2}(x,y)=f(y)-f(x). \label{2.28} \end{equation} Using notations \eqref{2.28}
we represent \eqref{2.26} in the form
\begin{equation}J_{1}(f)=-\int_{-1}^{0}\frac{\partial}{\partial{y}}P_{1}(x,y+x)\mu_{-}(y)dy
-\int_{-a-x}^{-1}\frac{\partial}{\partial{y}}P_{2}(x,y+x)\mu_{-}(y)dy, \label{2.29} \end{equation}Integrating by parts the integrals of \eqref{2.29}
 we deduce that
\begin{equation}J_{1}(f)=f^{\prime}(x)\gamma_{1}+\int_{-1}^{0}P_{1}(x,y+x)d\mu_{-}(y)
+\int_{-a-x}^{-1}P_{2}(x,y+x)d\mu_{-}(y)+P_{2}(x,-a)\mu_{-}(-a), \label{2.30} \end{equation}
where $\gamma_{1}=k_{-}^{\prime}(-1).$ It follows from \eqref{1.3} that the integrals in \eqref{2.30} are absolutely convergent. Passing to the limit in \eqref{2.30},when $a{\to}+\infty$, and taking into account  \eqref{2.16} , \eqref{2.27} we have
 \begin{equation}
J_{1}(f)=\int_{-\infty}^{x}[f(y+x)-f(x)-y\frac{df(x)}{dx}1_{|y|{\leq}1}]d{\mu_{-}}(y)
+{\gamma}_{1}f^{\prime}(x). \label{2.31} \end{equation}In the same way it can be
proved that \begin{equation}
J_{2}(f)=\int_{x}^{\infty}[f(y+x)-f(x)-y\frac{df(x)}{dx}1_{|y|{\leq}1}]d{\mu_{+}}(y)
+{\gamma}_{2}f^{\prime}(x), \label{2.32} \end{equation}where
$\gamma_{2}=k_{+}^{\prime}(1).$ The relation \eqref{2.26} follows directly
from  \eqref{2.31} and \eqref{2.32}. Here $\Gamma=\gamma_{1}+\gamma_{2}.$ The lemma is
proved.
\end{proof}
\begin{Rk}\label{Remark 2.5.}
The operator $L_{0}f=\frac{d}{dx}f$ can be
represented in  form \eqref{1.9}, \eqref{1.10}, where \begin{equation}
S_{0}f=\int_{-\infty}^{\infty}p_{0}(x-y)f(y)dy, \label{2.33} \end{equation}
\begin{equation}p_{0}(x)=\frac{1}{2}\,\mathrm{sign}(x). \label{2.34} \end{equation}
\end{Rk}

From Lemmas \ref{Lemma 2.4.},  and Remark \ref{Remark 2.5.} we deduce the following
assertion.

\begin{Tm}\label{Theorem 2.5}
The infinitesimal generator  $L$
has a
convolution-type form \eqref{1.9},  \eqref{1.10}.
\end{Tm}
\begin{Ee}\label{Example 2.7}Let us consider the Poisson process.\end{Ee} In this case we have
\begin{equation}\lambda(z)=-(e^{iz}-1)\label{2.35}\end{equation}
According to \eqref{1.2} and \eqref{2.35} the relations
\begin{equation}\mu_{-}(x)=0,\,x<0;\mu_{+}(x)=\begin{cases}-1 &\text{if $0<x<1$;}\\
0 &\text{if $x{\geq}1$} \label {2.36}\end{cases}\end{equation}
are true. Using \eqref{2.8},\eqref{2.9} and \eqref{2.36} we obtain
\begin{equation}k_{-}(x)=0,\,x<0;k_{+}(x)=\begin{cases}1-x &\text{if $0<x<1$;}\\
0 &\text{if $x{\geq}1$} \label {2.37}\end{cases}\end{equation}
Hence the operator L for the Poisson process has the following convolution form:
\begin{equation}Lf=\frac{d}{dx}\int_{x}^{x+1}(1-y+x)f^{\prime}(y)dy.\label{2.38}\end{equation}
Formula \eqref{2.38} coincides with Levy-Ito formula:
\begin{equation}Lf=-f^{\prime}(x) +f(x+1)-f(x).\label{2.39}\end{equation}
\section{The Probability of the Levy process (type II) remaining within the given domain}
\label{sec 3}
%%%%%
1.  We remind, that the definition of the Levy processes and the definition of the type II are given
in the section 1.\\Let us denote by $\Delta$ the set of segments $[a_{k},b_{k}]$ such that\\
 $a_{1}<b_{1}<a_{2}<b_{2}<...<a_{n}<b_{n},\quad 1{\leq}k{\leq}n.$
In many theoretical and applied problems it is important to estimate the quantity
\begin{equation}p(t,\Delta)=P(X_{\tau}{\in}\Delta; 0{\leq}\tau{\leq}t), \label{3.1} \end{equation}
i.e. the probability that a sample of the process $X_{\tau}$ remains
inside $\Delta$
for $0{\leq}\tau{\leq}t$ (\emph{ruin problem}).\\
\textbf{Condition 3.1} \emph{Further we consider only the Levy processes of type II and assume, that $\Delta$ belongs to the support of $X_t,\,t>0$ (see section 6).}\\
According to Condition 3.1 $\Delta{\in}\textbf{R}$ if either condition 1) or condition 2)
of Theorem 6.3 is fulfilled. If condition 3) of Theorem 6.3 is fulfilled and $\gamma=0$,
then either $\Delta{\in}[0,\infty)$ or $\Delta{\in}(-\infty,0].$\\
We denote by $F_{0}(x,t)$  the distribution function of Levy process $X_t$,i.e.
   \begin{equation}F_{0}(x,t)=P(X_{t}{\leq}x).
 \label{3.2} \end{equation}We need the following statement (see \cite{Sato})
\begin{Tm}\label{Theorem 3.1} The distribution function $F_{0}(x,t)$
 is continuous with respect to $x$ if and only if the Levy process belongs to type
 $II$.\end{Tm}
We introduce the sequence of functions
\begin{equation}
F_{n+1}(x,t)=\int_{0}^{t}\int_{-\infty}^{\infty}F_{0}(x-\xi,t-\tau)V(\xi)d_{\xi}F_{n}(\xi,\tau)d\tau,
 \label{3.3} \end{equation} where the function $V(x)$ is defined by relations
 $V(x)=1$ when $x{\notin}\Delta$ and $V(x)=0$ when $x{\in}\Delta.$
 In the right side of \eqref{3.3} we use Stieltjes integration. It follows from \eqref{1.1} that
\begin{equation}\mu(z,t)=\mu(z,t-\tau)\mu(z,\tau).\label{3.4}\end{equation}
Due to \eqref{3.4} and convolution formula for Stieltjes-Fourier transform (\cite{Boch}, Ch.4)
the relation
 \begin{equation}
F_{0}(x,t)=\int_{-\infty}^{\infty}F_{0}(x-\xi,t-\tau)d_{\xi}F_{0}(\xi,\tau) \label{3.5} \end{equation} is true.
 Using \eqref{3.3} and \eqref{3.5} we have \begin{equation}
0{\leq}d_{x}F_{n}(x,t){\leq}t^{n}d_{x}F_{0}(x,t)/n!,\,if\,dx>0. \label{3.6} \end{equation}
Relation \eqref{3.6} implies that
\begin{equation}
0{\leq}F_{n}(x,t){\leq}t^{n}F_{0}(x,t)/n!. \label{3.7}\end{equation}
 Hence the
series \begin{equation}
F(x,t,u)=\sum_{n-0}^{\infty}(-1)^{n}u^{n}F_{n}(x,t)  \label{3.8} \end{equation}
converges. The probabilistic meaning of $F(x,t,u)$ is defined by the
relation (see \cite{Kac1}, Ch.4):
\begin{equation}E\{\mathrm{exp}[-u\int_{0}^{t}V(X_{\tau})d\tau],c_{1}<X_{t}<c_{2}\}=
F(c_2,t,u)-F(c_1,t,u). \label{3.9} \end{equation}The inequality $V(x){\geq}0$ and
relation \eqref{3.9} imply that the function $F(x,t,u)$ monotonically decreases with respect
 to the variable "u" and  monotonically increases with respect
 to the variable $"x"$. Hence, the following formula
  \begin{equation}
 0{\leq}F(x,t,u){\leq}F(x,t,0)=F_{0}(x,t) \label{3.10} \end{equation}
is true. In view of \eqref{3.2} and \eqref{3.10} the Laplace transform
\begin{equation}
\Psi(x,s,u)=\int_{0}^{\infty}e^{-st}F(x,t,u)dt,\quad
s>0. \label{3.11} \end{equation} has  meaning. According to \eqref{3.3} the function
$F(x,t,u)$ is the solution of the equation
\begin{equation}
F(x,t,u)+u\int_{0}^{t}\int_{-\infty}^{\infty}F_{0}(x-\xi,t-\tau)V(\xi)d_{\xi}F(\xi,\tau,u)d\tau=
F_{0}(x,t) \label{3.12} \end{equation}
Taking from both parts of \eqref{3.12} the Laplace transform  bearing in mind  \eqref{3.11} and using the convolution property (see \cite{Boch}, Ch.4) we obtain
\begin{equation}
\Psi(x,s,u)+u\int_{-\infty}^{\infty}\Psi_{0}(x-\xi,s)V(\xi)d_{\xi}\Psi(\xi,s,u)=
\Psi_{0}(x,s), \label{3.13} \end{equation}where
\begin{equation}
\Psi_{0}(x,s)=\int_{0}^{\infty}e^{-st}F_{0}(x,t)dt. \label{3.14} \end{equation}
It follows from \eqref{1.1} and \eqref{3.14} that
\begin{equation}\int_{-\infty}^{\infty}e^{ixp}d_{x}\Psi_{0}(x,s)=\frac{1}{s+\lambda(p)}.
\label{3.15}\end{equation}According to \eqref{3.13} and \eqref{3.14}
 we have \begin{equation}
\int_{-\infty}^{\infty}e^{ixp}[s+\lambda(p)+uV(x)]d_{x}\Psi(x,s,u)=1. \label{3.16} \end{equation}
 Now we introduce the function
\begin{equation}h(p)=\frac{1}{2\pi}\int_{\Delta}e^{-ixp}f(x)dx, \label{3.17} \end{equation}
where the function $f(x)$ belongs to $C_{\Delta}$. Multiplying both parts  of \eqref{3.16}
by $h(p)$ and integrating them with respect to $p\quad  (-\infty<p<\infty)$ we deduce
the equality
\begin{equation}
\int_{-\infty}^{\infty}\int_{-\infty}^{\infty}e^{ixp}[s+\lambda(p)]h(p)d_{x}\Psi(x,s,u)dp=f(0).
 \label{3.18} \end{equation} We have used the relations \begin{equation}
V(x)f(x)=0,\quad -\infty<x<\infty, \label{3.19} \end{equation}
\begin{equation}
\frac{1}{2\pi}\mathrm{lim}\int_{-N}^{N}\int_{\Delta}e^{-ixp}f(x)dxdp=f(0),\quad N{\to}\infty.
 \label{3.20} \end{equation}
Since the function $F(x,t,u)$ monotonically decreases with respect
to $"u"$
 this is also true for the function $\Psi(x,s,u)$ (see \eqref{3.11}). Hence there exist the limits
 \begin{equation}F_{\infty}(x,s)=\mathrm{\lim}F(x,s,u),\quad
 \Psi_{\infty}(x,s)=\mathrm{\lim}\Psi(x,s,u),\quad u{\to}\infty. \label{3.21} \end{equation}
 It follows from \eqref{3.8} that
 \begin{equation}p(t,\Delta)=P(X_{\tau}{\in}\Delta,\,0<\tau<t)=\int_{\Delta}d_{x}F_{\infty}(x,t).
 \label{3.22}\end{equation}
Hence we have
\begin{equation}
 \int_{0}^{\infty}e^{-st}p(t,\Delta)dt=\int_{\Delta}d_{x}\Psi_{\infty}(x,s). \label{3.23} \end{equation}
 Using   relations \eqref{1.2} and \eqref{1.6}
 we deduce that
 \begin{equation}\lambda(z)\int_{-\infty}^{\infty}e^{-i{z}\xi}f(\xi)d\xi=
 -\int_{-\infty}^{\infty}e^{-i{z}\xi}[Lf(\xi)]d\xi.\label{3.24}\end{equation}
 2. By $C_{\Delta}$
 we denote the set of functions $g(x)$ on $L^{2}(\Delta)$ such that
\begin{equation}g(a_{k})=g(b_{k})=g^{\prime}(a_{k})=g^{\prime}(b_{k})=0,\quad 1{\leq}k{\leq}n,\quad g^{\prime\prime}(x){\in}L^{p}(\Delta),\quad p>1.
 \label{3.25} \end{equation}We introduce the operator $P_{\Delta}$ by relation
$P_{\Delta}f(x)=f(x)$ if $x{\in}\Delta$ and $P_{\Delta}f(x)=0$ if $x{\notin}\Delta.$\\
\begin{Dn}\label{Definition 3.2}
The operator\begin{equation}
L_{\Delta}=P_{\Delta}LP_{\Delta}=\frac{d}{dx}S_{\Delta}\frac{d}{dx},
\,where\,S_{\Delta}=P_{\Delta}SP_{\Delta},  \label{3.26} \end{equation}
is called a  truncated generator.\end{Dn}
 Relations
 \eqref{3.18},\eqref{3.21}  and \eqref{3.24} imply the following assertion.
\begin{Tm}\label{Theorem 3.3} {Let $X_t$ be the  Levy process.If the corresponding distribution function $F(x,t)$ is continuous with respect to x
  then the relation  \begin{equation}
  \int_{\Delta}(sI-L_{\Delta})fd_{x}\Psi_{\infty}(x,s)=f(0) \label{3.27} \end{equation} is true.}
 \end{Tm}
\begin{Rk}\label{Remark 3.4}For the symmetric stable processes equality \eqref{3.27} was deduced by M.Kac
\cite{Kac1}. For the   Levy processes with continuous density the equality \eqref{3.27}
was obtained in the book
\cite{Sakh6}. Now we deduced equality \eqref{3.27} for the   Levy processes with continuous distribution.\end{Rk}
3. Relations \eqref{3.23} implies the following assertion
\begin{Pn}\label{Proposition 3.5}The function $\Psi_{\infty}(x,s)$ for all $s>0$
is monotonically increasing and continuous with respect to $x$.\end{Pn}
The behavior of  $\Psi_{\infty}(x,s)$ when $s=0$ we shall consider separately,
using the following Hengartner and Thedorescu result (\cite{HenTheo}, see \cite{Sato}
too):
\begin{Tm}\label{Theorem 3.6}  Let $X_t$ be a Levy process.Then for any finite interval K
the estimation
\begin{equation}P(X_t{\in}K)=0(t^{-1/2}) \quad as\quad t{\to}\infty \label{3.28}\end{equation}
is valid.\end{Tm}Hence, we have the assertion
\begin{Tm}\label{Theorem 3.7}Let $X_t$ be a Levy process.Then for any integer ${n>o}$
the estimation
\begin{equation}p(t,\Delta)=0(t^{-n/2}) \quad as\quad t{\to}\infty \label{3.29}\end{equation}
is valid.\end{Tm}
\begin{proof} Let the interval $K$ be such, that $\Delta{\in}K$. Then the inequality
\begin{equation}P(X_t{\in}\Delta){\leq}P(X_t{\in}K)\quad\label{3.30}\end{equation}
holds. According to Levy processes properties (independent and stationary increments)
and \eqref{3.30} the following inequality
\begin{equation}p(t,\Delta){\leq}P^{n}(X_{t/n}{\in}K)=0(t^{-n/2}) \quad as\quad t{\to}\infty\quad\label{3.31}\end{equation}is true.\end{proof}
It follows from \eqref{3.9} that
\begin{equation}
 \int_{0}^{\infty}e^{-st}p(t,\delta)dt=\int_{\delta}d_{x}\Psi_{\infty}(x,s), \label{3.32} \end{equation}
where $\delta$ is a set of segments  which belong to $\Delta$ and
\begin{equation}p(t,\delta)=P(X_{\tau}{\in}\delta; 0{\leq}\tau{\leq}t). \label{3.33} \end{equation} If $\delta=\Delta$, then formula \eqref{3.32} coincides with formula \eqref{3.23} We need the following partial case of \eqref{3.32}:
\begin{equation}
 \int_{0}^{\infty}p(t,\delta)dt=\int_{\delta}d_{x}\Psi_{\infty}(x,0), \label{3.34} \end{equation} According to \eqref{3.31} the integral in the left side of \eqref{3.34} exists. Let us prove the following statement.
 \begin{Pn}\label{Proposition 3.8}The function $\Psi_{\infty}(x,0)$
is strictly  monotonically increasing and continuous in the domain $\Delta$.  \end{Pn}
\emph{Proof.} In view of \eqref{3.34} and Condition 3.1 the function $\Psi_{\infty}(x,0)$ is strictly  monotonically increasing. To prove that  $\Psi_{\infty}(x,0)$ is continuous we introduce the function
\begin{equation}q(t,\delta)=P(X_t{\in}\delta).\label{3.35}\end{equation}It is obvious
that
\begin{equation}p(t,\delta){\leq}q(t,\delta).\label{3.36}\end{equation}We consider $\delta=[x_1,x_2]$.
As $F(x,t)$ is continuous function with respect to $x$ the relations
\begin{equation}\lim_{x_2{\to}x_1}q(t,\delta)=0,\,\lim_{x_2{\to}x_1}p(t,\delta)=0 \label{3.37}\end{equation} are valid. Formulas \eqref{3.31},
\eqref{3.34} and \eqref{3.37} imply that the function $\Psi_{\infty}(x,0)$ is continuous.
The proposition is proved.
\begin{Rk}Formulas \eqref{3.32} and \eqref{3.34} will play in the next sections an essential role.\end{Rk}
\section{Quasi-potential}
\label{section 4}
1. We remind that the domain $\Delta$ is the set of segments $[a_{k},b_{k}]$ , where\\
 $a_{1}<b_{1}<a_{2}<b_{2}<...<a_{n}<b_{n},\quad 1{\leq}k{\leq}n.$\\
We denote by $D_{\Delta}$ the space of the continuous functions $g(x)$ on the domain $\Delta$ such that\\
$g(a_k)=g(b_k)=0,\quad 1{\leq}k{\leq}n.$\\
The norm in $D_{\Delta}$ is defined by the relation $||f||=sup_{x{\in}\Delta}|f(x)|.$
\begin{Dn}\label{Definition 4.1.}
The operator $B$ with the definition domain
 $D_\Delta$  is called \emph{a quasi-potential} if the
following relations
\begin{equation} -L_{\Delta}Bf=f,\quad f{\in}D_{\Delta}  \label{4.1}\end{equation}
\begin{equation} -BL_{\Delta}g=g,\quad g{\in}C_{\Delta} \label{4.2} \end{equation}
are true.
 \end{Dn}
\begin{Rk}\label{Remark 4.2}
In a number of cases (see the next section) we need relations \eqref{4.1} and \eqref{4.2}.
In these cases we can use the quasi-potential $B$, which is often simpler
than the corresponding potential $Q$.\end{Rk}
\begin{Tm}\label{Theorem 4.3} Let the considered Levy process $X_t$ belong to the
type II.
 Then the corresponding quasi-potential
 $B$ is bounded in the space $D_{\Delta}$, has  the form
\begin{equation}Bf=\int_{\Delta}f(y)d_{y}\Phi(x,y),\,f{\in}D_{\Delta},\label{4.3}
\end{equation}
where the real-valued function $\Phi(x,y)$ is continuous with respect to
$x$ and $y$, monotonically increasing with respect to $y.$
\end{Tm}
Before proving Theorem 4.3, we investigate some properties of the operator $B$, which is defined by \eqref{4.3}. The operator $B$  maps the space $D_{\Delta}$ into himself. We remind the following definitions.
\begin{Dn}\label{Definition 4.4} The total variation of a complex-valued function g, defined on $\Delta$ is the quantity\\
$$V_{\Delta}(g)=\sup_{P}\sum_{i=0}^{n_P-1}|g(x_{i+1})-g(x_{i})|,$$\\
where the supremum is taken over the set  of all partitions $P=(x_0,x_1,...,x_{n_P})$
of the $\Delta$.\end{Dn}
\begin{Dn}\label{Definition 4.5}. A complex-valued function g on the  $\Delta$ is said to be of bounded variation (BV function) on the $\Delta$ if its total variation is finite.
\end{Dn} By $D^*_{\Delta}$
we denote the conjugate space to $D_{\Delta}$. It is well-known that the space $D^*_{\Delta}$ consists from functions $g(x)$ with a bounded total variation $V_{\Delta}(g)$. The norm in $D^*_{\Delta}$ is defined by the relation $||g||=V_{\Delta}(g)$, the functional in $D_{\Delta}$ is defined by the relation
\begin{equation}(f,g)_{\Delta}=\int_{\Delta}f(x)d\overline{g(x)},\,f{\in}D_{\Delta},\,
g{\in}D^*_{\Delta}.\label{4.4}\end{equation} Hence, the conjugate operator
$B^*$ maps the space $D^*_{\Delta}$ into himself and has the form
\begin{equation}B^*g=\int_{\Delta}\Phi(y,x)dg(y).\label{4.5}\end{equation}

Let us consider an arbitrary inner point $y_0$ from the domain $\Delta.$ We introduce the new  domain
$\Delta(y_0)=\Delta-y_0$. We denote the corresponding truncated generator by
 $L_{\Delta}(y_0)$, the corresponding quasi-potential by $B(y_0)$, the corresponding
kernel by $\Phi(x,y,y_0)$, the corresponding $\Psi$ function by
 $\Psi_{\infty}(x,y,y_0)$
If $y_{0}=0$, then  according to
\eqref{4.5} we have
  $\Phi(0,x,0)=B(0)^{\star}\sigma(x),$  where $\sigma(x)=-1/2,$ when $x<0$ and
$\sigma(x)=1/2,$ when $x>0.$   The next assertion follows directly from \eqref{3.27}.
\begin{La}\label{Lemma 4.6}If $y_0=0$ belongs to the inner  part of  $\Delta$  and $\Phi(0,x,0)=\Psi_{\infty}(x,0,0)$, then
\begin{equation}-(L_{\Delta}(0)f,B(0)^*\sigma)_{\Delta}=f(0).
\label{4.6} \end{equation}
\end{La}Now we shall  reduce the general case to the case $y_{0}=0$.
 We introduce the  operator
\begin {equation} Uf=g(x),\,g(x)=f(x-y_0),\,x{\in}\Delta,
 \label{4.7}
  \end{equation} which maps the space $D_{\Delta(y_0)}$ onto
$D_{\Delta}$.  Using formulas \eqref{2.1} and \eqref{2.2} we
deduce that
\begin{equation} L_{\Delta}=UL_{\Delta}(y_0)U^{-1}. \label{4.8} \end{equation}
Hence the equality
\begin{equation}
B=UB(y_0)U^{-1} \label{4.9} \end{equation} is valid. The last equality can be
written in the terms of the kernels \begin{equation}
\Phi(x,y)=\Phi(x-y_{0},y-y_{0},y_0). \label{4.10} \end{equation}
Relation \eqref{3.27} in the case $\Delta(y_0)$ takes the form
\begin{equation}\int_{\Delta(y_0)}-L_{\Delta(y_0)}fd_{x}\Psi_{\infty}(x,0,y_0)=f(0).
\label{4.11}\end{equation}
According to \eqref{4.11} Lemma 4.6 can be written in the following form:
\begin{La}\label{Lemma 4.7}If  $\Phi(0,x,y_0)=\Psi_{\infty}(x,0,y_0)$, then
\begin{equation}-(L_{\Delta}(y_0)f,B(y_0)^*\sigma)_{\Delta(y_0)}=f(0).
\label{4.12} \end{equation}\end{La}
In view of \eqref{4.8} and \eqref{4.9} equality \eqref{4.12} can be rewritten in the form
\begin{equation}(-L_{\Delta(0)}g,B(0)^*\sigma(x-y_0))_{\Delta}=g(y_0),\label{4.13}\end{equation}
where
\begin{equation}g(x)=f(x+y_0). \label{4.14}\end{equation}
Using \eqref{4.12}
 we define the kernel $\Phi(x,y)$ of the operator $B$  by the relation
\begin{equation}\Phi(x,y)=\Psi_{\infty}(y-x,0,x),\label{4.15}\end{equation}
According to Proposition 3.8 and \eqref{4.15} we have the assertion
\begin{Pn}\label{Proposition 4.8}
 The function $\Phi(x,y)$ is continuous with respect to
$x$ and $y$ and   monotonically increasing with respect to $y.$
\end{Pn}
Relation \eqref{4.15} implies the equality
\begin{equation}-BL_{\Delta}g=g,\quad g{\in}C_{\Delta}.\label{4.16}\end{equation}
2. \emph{Sectorial properties}\\ We shall need the following Pringsheim's result .
 \begin{Tm}\label{Theorem 4.9} (see \cite{Titch}, Ch.1)
 Let $f(t)$ be
non-increasing function over $(0,\infty)$ and integrable on any
finite interval $(0,\ell)$. If $f(t){\to}0$ when
$t{\to}\infty$, then for any positive $x$ we have
\begin{equation}\frac{1}{2}[f(x+0)+f(x-0)]=\frac{2}{\pi}\int_{+0}^{\infty}\mathrm{cos}xu
[\int_{0}^{\infty}f(t)\mathrm{cos}tudt]du, \label{4.17} \end{equation}
\begin{equation}\frac{1}{2}[f(x+0)+f(x-0)]=\frac{2}{\pi}\int_{0}^{\infty}\mathrm{sin}xu
[\int_{0}^{\infty}f(t)\mathrm{sin}tudt]du. \label{4.18} \end{equation}
\end{Tm}We choose such the functions  $k_{-}(x)$ and $k_{+}(x)$ that
\begin{equation}k_{-}(-a)=k_{+}(a)=0.\,a>o.\label{4.19}\end{equation}
In this case we have
\begin{equation}k_{-}(x)=\int_{-a}^{x}\mu_{-}(t)dt,\,-\infty{\leq}x<0,\label{4.20}
\end{equation}
\begin{equation}k_{+}(x)=-\int_{x}^{a}\mu_{+}(t)dt,\,0<x{\leq}+\infty.\label{4.21}
\end{equation}
Using the integration by parts and taking into account \eqref{2.7} and \eqref{2.15} we deduce the assertion.
\begin{Pn}\label{Proposition 4.10}
Let conditions \eqref{4.19}- \eqref{4.21}
  be fulfilled. Then the relation
\begin{equation}
\int_{-a}^{a}k(t)\mathrm{cos}xtdt=
-[k_{+}^{\prime}(a)-k_{-}^{\prime}(-a)]\frac{1-\mathrm{cos}xa}{x^{2}}+\int_{-a}^{a}
\frac{1-\mathrm{cos}xt}{x^{2}}d{\mu}(t), \label{4.22} \end{equation} is true.
Here $\mu(t)=\mu_{-}(t)$ when $t<0$ and $\mu(t)=\mu_{+}(t)$ when $t>0$.
\end{Pn}
The following relations
\begin{equation}k_{-}^{\prime}(-a)=\mu_{-}(-a){\geq}0,\,k_{+}^{\prime}(a)=\mu_{+}(a){\leq}0
\label{4.23}\end{equation} are valid.
It follows from \eqref{4.23} and Theorem 4.9
that the kernel $k(x)$  admits the
representation
\begin{equation} k(x)=\frac{1}{2\pi}\int_{-\infty}^{\infty}m(t)e^{ixt}dt,\, \label{4.24} \end{equation}where
\begin{equation}\mathrm{Re}[m(t)]=-[k_{+}^{\prime}(a)-k_{-}^{\prime}(-a)]\frac{1-\mathrm{cos}xa}{x^{2}}+\int_{-a}^{a}
\frac{1-\mathrm{cos}xt}{x^{2}}d{\mu}(t){\geq}0. \label{4.25} \end{equation}
Further we use the following notions.
\begin{Dn}\label{Definition 4.11}
The bounded operator $S_{\Delta}$ in the space $L^{2}(\Delta)$ is called
\emph{sectorial} if \begin{equation}
(S_{\Delta}f,f){\ne}0,\quad f{\ne}0  \label{4.26} \end{equation}
and \begin{equation}
-\frac{\pi}{2}\beta{\leq}\mathrm{arg}(S_{\Delta}f,f){\leq}\frac{\pi}{2}\beta, \quad 0<\beta{\leq}1.
 \label{4.27} \end{equation}
  \end{Dn}
 \begin{Dn}\label{Definition 4.12}
 The sectorial operator $S$ is called
 \emph{a strongly sectorial}
 if for some $\beta<1$  relation
\eqref{4.27}
is valid.
 \end{Dn}
We note that we use two definitions for scalar product\\ $(f,g)=\int_{\Delta}f(x)\overline{g(x)}dx$ and  $(f,g)_{\Delta}=\int_{\Delta}f(x)d\overline{g(x)}.$ \\
3.Let the equality  $\Delta=[-a,a]$ holds. Then
due to \eqref{4.24}  the relation
\begin{equation}
(S_{\Delta}f,f)=\int_{-\infty}^{\infty}m(u)|\int_{\Delta}f(t)e^{iut}dt|^{2}du \label{4.28} \end{equation}
is valid. If $0{\in}{\Delta}$, then the  Levy measure $\nu(\Delta)>0$.
In this case the  entire function $m_{1}(x)=\Re{m(x)}$ is equal to zero  only in
finite number of points on every finite interval.
In view of \eqref{4.28} we have $(S_{\Delta}f,f){\ne}0$ and $\Re(S_{\Delta}f,f)>0,$ when $||f||{\ne}0.$
Thus, we have obtained the assertion.
\begin{Pn}\label{Proposition 4.13}If the Levy process $X_t$ belongs to the type II and $0$ belongs to $\Delta$,
then
\begin{equation}
-\frac{\pi}{2}<\mathrm{arg}(S_{\Delta}f,f)<\frac{\pi}{2},
\quad f(t){\in}L^{2}(\Delta),\,f{\ne}0.
 \label{4.29} \end{equation}\end{Pn}
 \begin{Rk}\label {Remark 4.14}We have deduced Proposition 4.13 for the case , when $\Delta=[-a,a].$
 It is easy to repeat this proof for an arbitrary $\Delta.$\end{Rk}
 4. We denote by $G$ the numerical range of $S_{\Delta}$. Let us consider the point $\alpha$, which belongs to the boundary of $G$.There exists a sequence of functions $f_{n}{\in}L^{2}(\Delta),\,||f_{n}||=1$, which converges in weakly sense to $f$ and
  $(S_{\Delta}f_{n},f_{n}){\to}=\alpha.$ In view of \eqref{1.3} the operator $S_{\Delta}$ is compact. So,
 $(S_{\Delta}f_{n},f_{n}){\to}(S_{\Delta}f,f)=\alpha.$
 The last relation implies:
 \begin{Pn}\label{Proposition 4.15} Let relations \eqref{4.29} be true. If $\alpha$ belongs to the closed set
 $\overline{G}$ then either $\alpha=0$, or $\Re{\alpha}>0.$\end{Pn}
 Now we can prove the following important in our theory result.
\begin{Tm}\label{Theorem 4.16}The operator $S_{\Delta}$ is strongly sectorial.\end{Tm}
\begin{proof} The assertion of the theorem in the case $A>0$ follows directly from
\eqref{4.29} and $\Re(S_{\Delta}f,f){\geq}A/2,\,||f||=1.$ According to theorem Toeplitz-
Hausdorff (see \cite{GuRao})  the numerical range set $G$ is convex. Using this fact and Proposition 4.15  we  obtain the assertion of the theorem  in the case A=0.\end{proof}
\begin{Rk}\label{Remark 4.17}Theorem 4.16 under additional conditions was  obtained before in our works (see \cite{Sakh6}).\end{Rk}
We note that operators $S_{\Delta}$
plays an important role not only in Levy processes theory (see \cite{Sakh5} and
\cite{Sakh6}).\\
5. According to \eqref{4.16} we have
\begin{equation}-L_{\Delta}BL_{\Delta}g=L_{\Delta}g,\quad g{\in}C_{\Delta}.\label{4.30}\end{equation}Using relations \eqref{3.25} we can represent the operator $L_{\Delta}$ in the form
\begin{equation}L_{\Delta}g=S_{\Delta}g^{\prime\prime},\quad g{\in}C_{\Delta}.\label{4.31}
\end{equation}It follows from Theorem 4.16 and equality \eqref{4.31}, that the range
$L_{\Delta}$ is dense in $D_{\Delta}$. Hence, it follows from \eqref{4.30}
the statement.
\begin{La}\label{Lemma 4.18}If the kernel $\Phi(x,y)$ satisfies the relation \eqref{4.15},
then
\begin{equation}-L_{\Delta}Bf=f,\quad f{\in}D_{\Delta}.\label{4.32}\end{equation}\end{La}
Now we can prove the Theorem 4.3.\\
\textbf{Proof of Theorem 4.3.}\\
The function   $\Phi(x,y)$, defined by (4.15),  is bounded,continuous and monotonically increasing with respect to $y.$ Hence, the corresponding operator $B$ is bounded in the space $D_{\Delta}.$ It follows from (4.16) and (4.32) that the constructed operator $B$ is quasi-potential.The theorem is proved.\\
6.We need the following result.
\begin{Pn}\label{Proposition 4.19}
{The operator $B$ is strongly sectorial.}
\end{Pn}
%\emph{Proof.}
\begin{proof}
Let the function $g(x)$ satisfies conditions \eqref{3.25}. Then the relation
\begin{equation}
(-L_{\Delta}g,g)=(S_{\Delta}g^{\prime},g^{\prime})  \label{4.33} \end{equation}
holds. Equalities \eqref{4.1} and \eqref{4.33} imply that
\begin{equation}(f,Bf)=(S_{\Delta}g^{\prime},g^{\prime}) ,\quad g=Bf. \label{4.34} \end{equation}
Inequality \begin{equation}(Bf,f){\ne}0,\,if\,f{\ne}0 \label{4.35}\end{equation} follows
from relation \eqref{4.34}.
The operator  $S_{\Delta}$ is strongly sectorial.Hence, according to \eqref{4.34}
the operator $B$ is strongly sectorial too.
 \end{proof}
6. Now we shall find the relation between $\Psi_{\infty}(x,s)$ and $\Phi(0,x).$ \begin{Tm}\label{Theorem 4.20} {Let the considered Levy process belong to the type
II.
  Then in the space $D_{\Delta}$ there is one and only one function
  \begin{equation}
  \Psi(x,s)=(I+sB^{\star})^{-1}\Phi(0,x)  , \label{4.36} \end{equation}
  which satisfies relation $\eqref{3.27}$.}
  \end{Tm}

  %\emph{Proof.}
  \begin{proof}
  In view of \eqref{4.2} we have   \begin{equation}
  -BL_{\Delta}f=f,\quad f{\in}C_{\Delta}.  \label{4.37} \end{equation}
  Relations \eqref{4.36} and \eqref{4.37}
  imply that \begin{equation}
((sI-L_{\Delta})f,\psi(x,s))_{\Delta}=-((I+sB)L_{\Delta}f,\psi)_{\Delta}=-(L_{\Delta}f,\Phi(0,x))_{\Delta}. \label{4.38} \end{equation}
Since $\Phi(0,x)=B^{\star}\sigma(x),$
 then according to \eqref{4.36} and \eqref{4.38} relation \eqref{3.27} is
true.\\
 Let us suppose that in $L(\Delta)$ there is another function
$\Psi_{1}(x,s)$ satisfying \eqref{3.27}. Then the equality
\begin{equation} ((sI-L_{\Delta})f,\phi(x,s))_{\Delta}=0,\quad
\phi=\Psi-\Psi_{1}  \label{4.39} \end{equation} is valid. We write  relation
\eqref{4.39} in the form
\begin{equation}
(L_{\Delta}f,(I+sB^{\star})\phi)_{\Delta}=0. \label{4.40} \end{equation}
 The range of $L_{\Delta}$ is dense in $D_{\Delta}.$ Hence
in view of \eqref{4.40} we have
$\phi=0.$ The theorem is proved.
\end{proof}
It follows from relations \eqref{3.27} and \eqref{4.36} that
\begin{equation}\Psi_{\infty}(x,s)=(I+sB^{\star})^{-1}\Phi(0,x)  , \label{4.41} \end{equation}
\begin{Rk}We stress an important fact: Operators $B$ and $B^{\star}$ are bounded
in the spaces $D_{\Delta}$ and $D_{\Delta}^{\star}$ respectively.\end{Rk}
7.\begin{Dn}The Levy process $X_t$ with the triplet $(A,\gamma,\nu)$ is symmetric
 (see \cite{Sato}) if $\gamma=0$ and the function $\nu(x)$ is odd.\end{Dn}
According to \eqref{2.2}, \eqref{2.3} and \eqref{2.12}, \eqref{2.13} in the symmetric
case the following equality $k_(-x)=k+(x),\,x>0$ is valid.
Taking into account \eqref{3.26} we obtain the assertion:
\begin{Pn}Let $X_t$ be the symmetric Levy process. Then the corresponding operators
$L_{\Delta}$ and $B$ are self-adjoint in the Hilbert space $L^{2}(\Delta)$
with the inner product
\begin{equation}(f,g) =\int_{\Delta}f(x)\overline{g(x)}dx.\label{4.42}\end{equation}
\end{Pn}
So, in the symmetric case the operator B is strongly sectorial and self-adjoint.
Hence we have
\begin{Cy}Let $X_t$ be the symmetric Levy process. The spectrum of the corresponding operator
$B$ is real and non-negative.\end{Cy}
\section{Long time behavior}\label{section 5}

1.We apply the following Krein-Rutman theorem (\cite{KR},section 6):\\
\begin{Tm}\label{Theorem 5.1} {If a linear compact operator $B$ leaving
invariant a cone $K$,  has a point of the spectrum different from
zero,  then it has a positive eigenvalue $\lambda_{1}$ not less in
modulus than any other eigenvalues $\lambda_{k},\quad (k>1)$.
 To this eigenvalue $\lambda_{1}$ corresponds
at least one eigenvector $g_{1}{\in}K, (Bg_{1}=\lambda_{1}g_{1})$ of
the operator $B$ and at least one eigenvector $h_{1}{\in}K^{\star},
(B^{\star}h_{1}=\lambda_{1}h_{1})$ of the operator $B^{\star}$.}
\end{Tm}
We remark that in our case the cone $K$ consists of non-negative
continuous real functions $g(x){\in}D_{\Delta}$ and the cone $K^{\star}$ consists of monotonically increasing bounded functions $h(x){\in}D_{\Delta}^{\star}.$
In this section we investigate the asymptotic
 behavior
   $p(t,\Delta)$ when $t{\to}\infty$.\\
 The spectrum
$(\lambda_{k},\, k>1)$ of the operator $B$ is situated
 in the sector
\begin{equation}
-\frac{\pi}{2}\beta{\leq}\mathrm{arg}z{\leq}\frac{\pi}{2}\beta,\quad 0{\leq}\beta<1,
\quad |z|{\leq}\lambda_{1}.
 \label{5.1} \end{equation}
We introduce the domain $D_{\epsilon}$:
\begin{equation}
-\frac{\pi}{2}(\beta
+\epsilon){\leq}\mathrm{arg}z{\leq}\frac{\pi}{2}(\beta
+\epsilon),\quad |z-(1/2)\lambda_{1}|<(1/2)(\lambda_{1}
-r), \label{5.2} \end{equation} where $0<\epsilon<1-\beta),\quad 0<r<\lambda_{1}$.
If $z$ belongs to the domain $D_{\epsilon}$ then the relation
\begin{equation} \mathrm{Re}(1/z)>1/\lambda_{1} \label{5.3} \end{equation}
holds. Indeed, relation \eqref{5.3} is equivalent to inequality
$$(x-\lambda_{1}/2)^{2}+y^{2}<\lambda_{1}^{2}/4,\,x=\Re{z},\,y=\Im{z},\,i.e.\,
|z-\lambda_{1}/2|<\lambda_{1}/2.$$ We take so small $r$ that the
 the circle
\begin{equation} |z-(1/2)\lambda_{1}|=(1/2)(\lambda_{1}-r).\label{5.4}\end{equation}
has the points $z_{1}$ and  $z_{2}= \overline{z_{1}}$ of the intersections with
half-lines $arg{z}=\frac{\pi}{2}(\beta+\epsilon)$ and $arg{z}=-\frac{\pi}{2}(\beta+\epsilon)$. We denote  the boundary of
domain $D_{\epsilon}$
 by $\Gamma_{\epsilon}$. We stress that $\Gamma_{\epsilon}$ contains only part of circle \eqref{5.4}, which is situated between the points $z_1$ and $z_2$. Without loss of generality we may assume
 that the points of spectrum
 $\lambda_{k}{\ne}0$ do not belong to $\Gamma_{\epsilon}$.
 Now we formulate the main result of this section.
\begin{Tm}\label{Theorem 5.2}
Let Levy process $X_t$ have  type  $II$ and let the corresponding
quasi-potential $B$ satisfy the following conditions:\\
I. Operator $B$ is compact in the Banach space $D_{\Delta}$.
II. There exists  such constant $M>0$ that
\begin{equation}|((B-zI)^{-1}1,\Phi(0,x))_{\Delta}|{\leq}M/|z|.\,z{\in}\Gamma_{\ve}|.
\label{5.5}\end{equation}
Then we have\\
1).The eigenfunction $g_1(x)$ of the operator $B$ is continuous and $g_1(x)>0$, where
$x$ are the   inner points of $\Delta.$\\
2).The eigenfunction $h_1(x)$ of the operator $B^{\star}$ is absolutely continuous and strictly
monotonic.\\
3).The asymptotic equality \begin{equation}
 p(t,\Delta)=e^{-t/\lambda_{1}}[q(t)+o(1)],\quad t{\to}+\infty
  \label{5.6} \end{equation}
  {is true. The function $q(t)$ has the form}
 \begin{equation}
 q(t)=c_{1}+\sum_{k=2}^{m}c_{k}e^{it\nu_{k}}{\geq}0,\,c_{1}>0, \label{5.7} \end{equation}
{where $\nu_{k}$ are real.}
\end{Tm}

\begin{proof}The operator $B$  has  a point of the
 spectrum different from zero (see Theorem 6.6).\\
Properties 1) and 2) follows directly from Krein-Rutman Theorem 5.1. Hence we have
\begin{equation}(g_1,h_1)_{\Delta} =\int_{\Delta}g_1(x)dh_1(x)>0.\label{5.8}\end{equation} Using \eqref{3.23}  we obtain   the equality
\begin{equation} p(t,\Delta)=\frac{1}{2\pi}\int_{-\infty}^{\infty}(e^{iyt},\Psi_{\infty}(x,iy))_{\Delta}dy,
\,t>0.\label{5.9}\end{equation} Changing the variable $z=i/y$ and taking into account
\eqref{4.36} we rewrite \eqref{5.9} in the form
\begin{equation} p(t,\Delta)=\frac{1}{2i\pi}\int_{-i\infty}^{i\infty}(e^{-t/z},(zI-B^{\star})^{-1}
\Phi(0,x))_{\Delta}\frac{dz}{z},
\,t>0.\label{5.10}\end{equation}
As the operator $B$ is compact, only a finite number of
eigenvalues $\lambda_{k}, \quad 1<k{\leq}m$ of this operator does
not belong to the domain $D_{\epsilon}$.
 We deduce from formula \eqref{5.10}  the relation
 \begin{equation}
 p(t,\Delta)=
 \sum_{k=1}^{m}\sum_{j=0}^{n_{k}-1}e^{-t/\lambda_{k}}t^{j}c_{k,j}+J, \label{5.11} \end{equation}
where $n_{k}$ is the index of the eigenvalue
$\lambda_{k}$,
\begin{equation}
J=-\frac{1}{2i\pi}\int_{\Gamma_{\ve}}\frac{1}{z}e^{-t/z}(1,(B^{\star} -zI)^{-1}\Phi(0,x))_{\Delta}dz.
 \label{5.12} \end{equation}
 We remind that the \emph{index} of the  eigenvalue
$\lambda_{k}$ is defined as the dimension of the largest Jordan block associated to that eigenvalue.
We note that \begin{equation}n_{1}=1. \label{5.13} \end{equation}
Indeed, if $n_{1}>1$ then there exists such a function $f_{1}$ that
\begin{equation}
Bf_{1}=\lambda_{1}f_{1}+g_{1}. \label{5.14} \end{equation}
In this case the relations
\begin{equation}
(Bf_{1},h_{1})_{\Delta}=\lambda_{1}(f_{1},h_{1})_{\Delta}+(g_{1},h_{1})_{\Delta}
=\lambda_{1}(f_{1},h_{1})_{\Delta}
 \label{5.15} \end{equation}are true. Hence $(g_{1},h_{1})_{\Delta}=0.$ The last relation contradicts
 \eqref{5.8}. It proves equality \eqref{5.13}.\\
 Relation \eqref{4.13} implies  that
 \begin{equation}
 \Phi(0,x){\in}D^{\star}_{\Delta}. \label{5.16} \end{equation}
Among the numbers $\lambda_{k}$ we choose the ones for which
$\mathrm{Re}(1/{\lambda}_{k}),\quad (1{\leq}k{\leq}m)$ has the
smallest value $\delta$. Among the obtained numbers we choose
$\mu_{k},\quad (1{\leq}k{\leq}\ell)$ the indexes $n_{k}$ of which
have the largest value $n$. We deduce from \eqref{5.11}, \eqref{5.12} and \eqref{5.5}  that
\begin{equation}
 p(t,\Delta)=e^{-t{\delta}}t^{n}
 [\sum_{k=1}^{\ell}e^{-t/{\mu}_{k}}c_{k}+o(1)],\quad t{\to}\infty. \label{5.17} \end{equation}
 We note that the function
 \begin{equation} Q(t)= \sum_{k=1}^{\ell}e^{it\mathrm{Im}(\mu_{k}^{-1})}c_{k}  \label{5.18} \end{equation}
 is almost periodic (see \cite{LEV}). Hence in view of \eqref{5.17} and the inequality\\
  $p(t,\Delta)>0,\quad (t{\geq}0)$ the following relation
 \begin{equation} Q(t){\geq}0,\quad -\infty<t<\infty  \label{5.19} \end{equation}
 is valid.\\
  First we assume that at least one of the inequalities
 \begin{equation}\delta<{\lambda}_{1}^{-1},\quad n>1 \label{5.20} \end{equation}
 is true. Using  \eqref{5.20} and the inequality
  \begin{equation} \lambda_{1}{\geq}|\lambda_{k}|,\quad k=2,3,... \label{5.21} \end{equation}we have
\begin{equation}\mathrm{Im}\mu_{j}^{-1}{\ne}0,\quad 1{\leq}j{\leq}\ell. \label{5.22} \end{equation}
It follows from \eqref{5.18} that
\begin{equation}
c_{j}=\mathrm{lim}\frac{1}{2T}\int_{-T}^{T}Q(t)e^{-it(\mathrm{Im}{\mu}_{j}^{-1})}dt,\quad
T{\to}\infty. \label{5.23} \end{equation}In view of \eqref{5.19}  the relations
\begin{equation}
|c_{j}|{\leq}\mathrm{lim}\frac{1}{2T}\int_{-T}^{T}Q(t)dt=0,\quad
T{\to}\infty, \label{5.24} \end{equation} are valid, i.e. $c_{j}=0,\quad
1{\leq}j{\leq}\ell.$ This means that  relations \eqref{5.18} are not true.
Hence the equalities
 \begin{equation}\delta={\lambda}_{1}^{-1},\quad n=1 \label{5.25} \end{equation}
 are true. From  \eqref{5.25} we get the asymptotic equality
 \begin{equation}
 p(t,\Delta)=e^{-t/\lambda_{1}}[q(t)+o(1)]\quad t{\to}\infty, \label{5.26} \end{equation}
where the function $q(t)$ is defined by relation \eqref{5.7} and
 \begin{equation} c_{k}=g_{k}(0)\int_{\Delta}d\overline{h_{k}(x)},
 \quad \nu_{k}=\mathrm{Im}(\mu^{-1}),\,k>1, \label{5.27} \end{equation}
 \begin{equation} c_{1}=g_{1}(0)\int_{\Delta}dh_{1}(x)>0.\label{5.28}\end{equation}
 Here $g_{k}(x)$ are the eigenfunctions of the operator $B$ corresponding to
  the eigenvalues $\lambda_{k}$, and $h_{k}(x)$ are the eigenfunctions of the operator $B^{\star}$
   corresponding to
  the eigenvalues $\overline{\lambda_{k}}.$ The following conditions are fulfilled
  \begin{equation}
 (g_{k},h_{k})=\int_{\Delta}g_{k}(x)d\overline{h_{k}(x)}=1, \label{5.29} \end{equation}
\begin{equation}
 (g_{k},h_{\ell})=\int_{\Delta}g_{k}(x)d\overline{h_{\ell}(x)}=0,\quad k{\ne}\ell. \label{5.30} \end{equation}Using the almost periodicity of the function $q(t)$ we deduce from \eqref{5.17} the inequality $q(t){\geq}0$. The theorem is proved.
\end{proof}

 2. Now we shall consider the important case when
\begin{equation}\mathrm{rank}\lambda_{1}=1. \label{5.31} \end{equation}
We remain that \emph{rank} of an eigenvalue is defined as the number of linearly
independent eigenvectors with that eigenvalue, i.e. rank  of an eigenvalue coincides
with the geometric multiplicity of this eigenvalue.
\begin{Tm}\label{Theorem 5.3}
{Let the conditions of Theorem \ref{Theorem 5.2} be
fulfilled. In the case \eqref{5.31} the following relation}
\begin{equation}
 p(t,\Delta)=e^{-t/\lambda_{1}}[c_{1}+o(1)],\,c_1>0,\,\quad t{\to}+\infty
  \label{5.32}
  \end{equation}
 {is true.}
\end{Tm}

\emph{Proof}
In view of \eqref{5.6}
 we have
\begin{equation}
\mathrm{lim}\frac{1}{T}\int_{0}^{T}q(t)dt{\geq}
|\mathrm{lim}\frac{1}{T}\int_{0}^{T}q(t)e^{-it(\mathrm{Im}{\mu}_{j}^{-1})}dt|,\quad
T{\to}\infty, \label{5.33} \end{equation}i.e.

\begin{equation}
g_{1}(0)\int_{\Delta}dh_{1}(x){\geq}|g_{j}(0)\int_{\Delta}d\overline{h_{j}(x)}|. \label{5.34} \end{equation}
In the same way we can prove that
\begin{equation}
g_{1}(x_{0})\int_{\Delta_{1}}dh_{1}(x){\geq}|g_{j}(x_{0})\int_{\Delta_{1}}d\overline{h_{j}(x)}|,
\label{5.35}\end{equation}where

\begin{equation}x_{0}{\in}\Delta_{1}{\in}\Delta. \label{5.36} \end{equation} It follows from \eqref{5.35} that

 \begin{equation}g_{1}(x_{0})h_{1}(x){\geq}|g_{j}(x_{0})\overline{h_{j}(x)}|.\label{5.37}
  \end{equation}
 We introduce the normalization conditions
\begin{equation}g_{1}(x_{0})= g_{j}(x_{0}),\,h_{1}(x_{0})=h_{j}(x_{0})=0.\label{5.38}
\end{equation}

 Due to \eqref{5.35} and \eqref{5.37} the inequalities
\begin{equation}
\int_{\Delta_{1}}dh_{1}(x){\geq}|\int_{\Delta_{1}}dh_{j}(x)|.
\label{5.39}
 \end{equation}
 \begin{equation}
h_{1}(x){\geq}|h_{j}(x)| \label{5.40} \end{equation} are true. The equality sign
in \eqref{5.37} and \eqref{5.38} can be only if
\begin{equation}
h_{j}(x)=|h_{j}(x)|e^{i\alpha}. \label{5.41} \end{equation} It is possible only in the case when $j=1$. Hence
there exists such a
point $x_{1}$ that
\begin{equation}
h_{1}(x_{1})>|h_{j}(x_{1})| \label{5.42} \end{equation} Thus we have
\begin{equation}1=\int_{\Delta_{1}}g_{1}(x)dh_{1}(x)>
\int_{\Delta_{1}}g_{j}(x)d\overline{h_{j}(x)}=1,
 \label{5.43} \end{equation}where $x_{1}{\in}\Delta_{1}.$
 The received contradiction \eqref{5.43} means that $j=1.$ Now the assertion
of the theorem follows directly from \eqref{5.2} .
\begin{cor}\label{Corollary 5.2.}
{Let conditions of Theorem \ref{Theorem 5.2} be
fulfilled. If $\mathrm{rank}\lambda_{1}=1$ and
$x_{0}{\in}\Delta_{1}{\in}\Delta$ then
 the asymptotic equality }
 \begin{equation}
 p(x_{0},\Delta_{1},t,\Delta)=e^{-t/\lambda_{1}}g_{1}(x_{0})\int_{\Delta_{1}}h_{1}(x)dx[1+o(1)],\quad t{\to}+\infty
  \label{5.44} \end{equation}
  { is true.}
  \end{cor}
 According to Theorem \ref{Theorem 5.2} and the relation $0<Re{(1/\lambda_{k})}{\leq}1/\lambda_{1}$ the following assertion
  is true.
 \begin{cor}\label{Corollary 5.2}
 { Let the conditions of  Theorem \ref{Theorem 5.2}  be fulfilled.
 Then all the eigenvalues $\lambda_{j}$ of $B$
belong  to the disk}
 \begin{equation}
|z-(1/2)\lambda_{1}|{\leq}(1/2)\lambda_{1}. \label{5.45} \end{equation}
{All the eigenvalues $\lambda_{j}$ of $B$ which belong to the boundary of
 disc \eqref{5.45} have the indexes} $n_{j}=1.$
 \end{cor}
3. We consider again the Levy process $X_t$ of the type II. Now we do not suppose
that the quasi-potential operator $B$ is compact. As we proved before the operator $B$
is bounded in the Banach space $D_{\Delta}$.
\begin{Tm}\label{Theorem 5.4}Let the considered Levy process $X_t$ have the type II and let the spectrum of the corresponding operator $B$ belong to the intersection of the discs
 $|\lambda|{\leq}\rho$ and $|z-(1/2)\rho|{\leq}(1/2)\rho$. If the condition \eqref{5.5} is fulfilled, then we have
\begin{equation}p(t,\Delta)=O(e^{-t/(\rho+\ve)}),\,\,\ve>0,\,t{\to}\infty,\label{5.46}\end {equation}
\end{Tm}
4.Now we assume that the function $\Phi(x,y)$ is absolutely continuos with respect to$y$.
Hence, the corresponding operator $B$ has the form
\begin{equation}Bf=\int_{\Delta}f(y)K(x,y)dy,\quad K(x,y)=\frac{\partial}{\partial{y}}\Phi(x,y),\label{5.47}\end{equation}
where for all $x{\in}\Delta$ the inequality
\begin{equation}\int_{\Delta}|K(x,y)|dy<\infty \label{5.48}\end{equation}
holds.In addition we assume that for all $x{\in}\Delta$ the following relation
\begin{equation}\lim_{h{\to}0}\int_{\Delta}|K(x+h,y)-K(x,y)|dy=0 \label{5.49}\end{equation}
is fulfilled. Then the operator $B$ is compact in the space $D_{\Delta}$.
Here we use  the well-known  Radon's theorem \cite{RAD}:
\begin{Tm}\label{Theorem 5.5}The operator $B$ of the form \eqref{5.47} is compact in the space $D_{\Delta}$
if and only if the relations \eqref{5.48} and \eqref{5.49} are fulfilled.\end{Tm}
5.We can formulate the condition \eqref{5.5} of Theorem 5.2  in the terms of the kernel $\Phi(x,y)$.
We suppose that there exists such positive, monotonically decreasing function r(s) that
\begin{equation}|K(x,y)|{\leq}r(|x-y|),\,\int_{0}^{b}r^{2}(s)ds<\infty,\,b>0.\label{5.50}
\end{equation}It follows from \eqref{5.50}, that the operator $B$ is bounded in the Hilbert space $L^{2}(\Delta).$  The numerical range $W(B)$ of the operator $B$ in the Hilbert space $L^{2}(\Delta)$ is the set $W(B)=\{(Bf,f),\,||f||=1\}.$
The closure of the convex hull of $W(B)$ is situated in the sector \eqref{5.1}. Hence,
the estimation
\begin{equation}||(B-zI)^{-1}1||{\leq}M/|z|.\,z{\in}\Gamma_{\ve}
\label{5.51}\end{equation} is valid (see \cite{STO}).  It follows from \eqref{5.50}
and \eqref{5.51} that
\begin{equation}|((B-zI)^{-1}1,K(0,x))|{\leq}M/|z|.\,z{\in}\Gamma_{\ve}
\label{5.52}\end{equation} So, we have proved the following statement.
\begin{Tm}\label{Theorem 5.6} Let relations \eqref{5.49} and \eqref{5.50}
be fulfilled.Then the inequality \eqref{5.5} holds.\end{Tm}

\section{Appendix}\label{section 6}
1.\textbf{Support.}\\
Here we remind the following  definitions  (see \cite{Sato}).
\begin{Dn}\label{Definition 6.1}For any measure $\rho$ on $\mathbf{R}$ its support
$S(\rho)$ is defined to be the set of $x{\in}\mathbf{R}$ such that $\rho(G)>0$ for any
interval $G$ containing $x$.\end{Dn}
The support $S(\rho)$ is closed set.
\begin{Dn}\label{Definition 6.2}For any random variable $X$ on $\mathbf{R}$ the  support
of the corresponding distribution $F(x)$ is called the support of $X$ and denoted by
$S(X)$.\end{Dn}
Now we can formulate the well-known H.Tuckers theorem (\cite{Tuc}).
\begin{Tm}\label{Theorem 6.3}Let $X_t$ be a Levy process on \textbf{R} with triplet $(A,\nu,\gamma)$.\\
1) If either $A>0$ or
\begin{equation}\int_{-\infty}^{\infty}|x|d\nu=\infty, \label{6.1}\end{equation}
then $S(X_t)=\textbf{R}.$\\
2)  If $0{\in}S(\nu),\,S(\nu){\bigcap}(0,\infty){\ne}0,\,S(\nu){\bigcap}(-\infty,0){\ne}0,$ then $S(X_t)=\textbf{R}.$\\
3) Suppose that $A=0,\,0{\in}S(\nu)$ and
\begin{equation}\int_{-\infty}^{\infty}|x|d\nu<\infty. \label{6.2}\end{equation}
If $S(\nu){\in}[0,\infty),$ then $S(X_t)=[t\gamma,\infty).$\\
If $S(\nu){\in}(-\infty,0],$ then $S(X_t)=(-\infty,t\gamma].$
\end{Tm}
2.\textbf{Continuity.}\\
Let us consider the Levy process $X_t$ of the type II.
In this case the corresponding distribution $F(x,t)$ is continuous with respect to $x$
(see Theorem 1.3). We use the representation  (see\eqref{1.1}):
\begin{equation}\int_{-\infty}^{\infty}e^{ixz}dF(x,t)=\mu(z,t).\label{6.3}\end{equation}Hence,
the equality
\begin{equation}\int_{0}^{x}[F(y,t)-F(-y,t)]dy=\int_{\infty}^{\infty}\mu(\xi,t)\frac{1-cosx\xi}{\xi^{2}}d\xi.
\label{6.4}\end{equation} holds (\cite{Zy},Ch.XVI.) Now we shall prove the following fact
\begin{Tm}\label{Theorem 6.4} Let $X_t$ be a Levy process of type II. Then the corresponding
distribution function $F(x,t)$ is jointly continuous with respect to $t>0$
and $x{\in}\textbf{R}$.\end{Tm}
\begin{proof}Using \eqref{6.4} and the dominated convergence theorem we deduce that
the left side of \eqref{6.4} is continuous with respect to $t>0$
and $x{\in}\textbf{R}$. Hence, the function
\begin{equation} G(x,t)=F(x,t)-F(-x,t) \label{6.5}\end{equation}
is continuous with respect to $t>0$. Suppose , that $X_t$ satisfies condition 3)
of Theorem 6.3 and $\gamma=0$. In this case we have either $F(-x,t)=0,\,x>0$
or $F(x,t)=0,\,x>0$. As the function  $ G(x,t)$ is continuous with respect to $t>0$
the function $F(x,t)$ is continuous with respect to $t>0$ too.
Now let us consider an arbitrary Levy process $X_t$ of type II.
We can represent  $X_t$ in the following form
\begin{equation} X_t= X_{t}^{(1)}+X_{t}^{(2)},\label{6.6}\end{equation}
where Levy processes  $X_{t}^{(1)}$ and $X_{t}^{(2)}$ are independent and  the process
$X_{t}^{(2)}$ satisfies  condition 3)
of Theorem 6.3. The processes $X_{t}^{(1)}$ and $X_{t}^{(2)}$ have the following
Levy triplets $(A,\gamma,\nu_{1})$ and $(0,0,\nu_{2}$,
where
\begin{equation}\nu(x)=\nu_{1}(x)+\nu_{2}(x),\,\int_{-\infty}^{\infty}d\nu_{2}(x)=\infty.
\label{6.7}\end{equation} The distribution functions of $X_t$, $X_{t}^{(1)}$ and
$X_{t}^{(2)}$ are denoted by $F(x,t)$,  $F_{1}(x,t)$ and  $F_{2}(x,t)$ respectively.
The convolution formula for Stieltjes-Fourier transform (\cite{Boch}, Ch.4)
gives
 \begin{equation}
F(x,t)=\int_{-\infty}^{\infty}F_{2}(x-\xi,t-\tau)d_{\xi}F_{1}(\xi,\tau) \label{6.8}
\end{equation}Since the process  $X_{t}^{(2)}$ is of type II and satisfies the condition
3) of Theorem 6.3, the distribution  $F_{2}(x,t)$ is continuous with respect to $x$ and
with respect to $t$. The assertion of the theorem follows from this fact and relation
\eqref{6.8}
\end{proof}
\begin{Ee}\label{Example 6.5} We consider the case when
\begin{equation}\int_{0}^{\infty}d\nu(x)=\infty.\label{6.9}\end{equation}
\end{Ee}We shall explain the method of constructing
Levy measures $\nu_{1}(x)$ and $\nu_{2}(x)$. We take
\begin{equation}\nu_{1}(x)=\nu(x)-\nu_{2}(x)
\label{6.10}\end{equation}
\begin{equation}\nu_{2}(x)=0\,(x<0);d\nu_{2}(x)=xd\nu(x)\,(0{\leq}x{\leq}1);
d\nu_{2}(x)=d\nu(x)\,(x>1).\label{6.11}\end{equation}
It is easy to see that the corresponding Levy processes $X_{t}^{(1)}$ and $X_{t}^{(2)}$
satisfy \eqref{6.6}, are independent and  the process
$X_{t}^{(2)}$ satisfies  condition 3).\\
3. \textbf{Non-zero points of the spectrum.}\\
\begin{Tm}\label{Theorem 6.6}
Let Levy process $X_t$ have  type  $II$ and let the corresponding
quasi-potential $B$
 is compact in the Banach space $D_{\Delta}$.
Then the operator $B$ has a point of the spectrum different from zero.
\end{Tm}
\begin{proof}Let  a continuous, non-negative function $u(x)$ in the domain
$\Delta$ be such that $||u(x)||=1$ and
\begin{equation}u(x)=0,\,x{\in}[a_k,a_k+{\ve}_k][b_k-{\ve}_k,b_k],\,0<{\ve}_k<(b_k-a_k)/2,\,
1{\leq}k{\leq}n.\label{6.12}\end{equation}
The function $v(x)=Bu(x)$ will be continuous in the domain $\Delta$ and
$v(x)>0$, where  $x$ are inner points of $\Delta$. Hence there exists such $c>0$ that
\begin{equation}v(x){\leq}cu(x).\label{6.13}\end{equation}
The assertion of the theorem follows directly from \eqref{6.13}
and Krein-Rutman result (see \cite{KR}, theorem 6.2).
\end{proof}
4.\textbf{Weakly singular operators}\\
The operator
\begin{equation}Bf=\int_{\Delta}K(x,y)f(y)dy \label{6.14}\end{equation} is weakly singular if the
kernel $K(x,y)$ satisfies the inequality
\begin{equation}|K(x,y)|{\leq}M|x-y|^{\beta},\,0<\beta<1,\,M=constant.\label{6.15}
\end{equation}
\begin{Pn}\label{Proposition 6.7}The weakly singular  operator $B$ is bounded and compact in the spaces
$L^{p}(\Delta),\,1{\leq}p{\leq}\infty.$\end{Pn}
\begin{proof} The   operator $B$ is bounded  in the spaces
$L^{p}(\Delta),\,1{\leq}p{\leq}\infty$ (see \cite{Sakh5},p.24).
A sufficiently high iterate of the kernel $K(x,y)$ is bounded.
Hence the  operator $B$ is  compact in the spaces
$L^{p}(\Delta),\,1{\leq}p{\leq}\infty.$
\end{proof}
In particular,the weakly singular  operator $B$ is bounded and compact in the Hilbert
space $L^{2}(\Delta).$  Repeating the arguments of the Theorem 5.2 in the case  $L^{2}(\Delta),$ we obtain the assertion.
\begin{Tm}\label{Theorem 6.8}
Let Levy process $X_t$ have  type  $II$ and let the corresponding
quasi-potential $B$ is weakly singular.
Then we have\\
1)\begin{equation} K(x,y)=\frac{\p}{\p{y}}\Phi(x,y).\label{6.16}\end{equation}
2) The eigenfunction $g_1(x)$ of the operator $B$ is continuous and $g_1(x)>0$, where
$x$ are the   inner points of $\Delta.$\\
3)The eigenfunction $H_1(x)=h^{\prime}_{1}(x){\in}L^{2}(\Delta)$ is almost everywhere positive.
4)The equalities \eqref{5.6} and \eqref{5.7} are valid.\end{Tm}
\begin{proof} Now we consider the operator $B$ in the Hilbert space $L^{2}(\Delta)$.  We note that the operator $B$ has the same eigenvalues $\lambda_j$
and the same eigenvectors $g_j$ in the spaces $L^{2}(\Delta)$ and $D_{\Delta}$.
The operator $B^{\star}$ has the same eigenvalues $\overline{\lambda_j}$
in the spaces $L^{2}(\Delta)$ and $D_{\Delta}^{\star}$. The corresponding eigenvectors
$H_j(x)$ and  $h_j(x)$ are connected by the relation
\begin{equation}H{_j}(x)=h^{\prime}_{j}(x).\label{6.17}\end{equation}
We denote by $W(B)$ the numerical range of $B$.The closure of the convex hull of
 $W(B)$ is situated in the sector \eqref{5.2}.Hence the estimation
\begin{equation}||(B-zI)^{-1}||{\leq}M/|z|,\,z{\in}\Gamma_{\ve} \label{6.18}\end{equation} is valid \cite{STO}. The contour $\Gamma_{\ve}$ is defined in the section 5. If $\beta<1/2$ then $K(0,x){\in}L^{2}(\Delta)$. According to \eqref{6.18} we have
\begin{equation}|((B-zI)^{-1}1,K(0,x))|{\leq}M/|z|,\,z{\in}\Gamma_{\ve} \label{6.19}\end{equation}All conditions of Theorem 5.2 are fulfilled. So, in the case
$\beta<1/2$ the theorem is proved.\\
Let us consider the case when $1/2{\leq}\beta1/2<1.$\\
In this case we have
\begin{equation}(B^{\star})^{m}K(0,x){\in}L(\Delta),\, m=0,1,2,...\label{6.20}\end{equation}
\end{proof}

%%%%%%%

\begin{thebibliography}{1}
%\thispagestyle{empty}\markboth{Bibliography}{Bibliography}
%\addcontentsline{toc}{chapter}{Bibliography}
%%1. \textbf
\bibitem{BaDo}
{ G. Baxter and M. D. Donsker,} \textit {On the Distribution of the
Supremum Functional for Processes with Stationary Independent
Increments.} Trans. Amer. Math. Soc. \textbf {8} (1957), 73--87.
%2. \textbf
\bibitem{Bert}
{ J. Bertoin,} \textit {Levy Processes.} University Press, Cambridge, 1996.
%2. \textbf
\bibitem{Boch}
{Dr.S.Bochner.}Vorlesungen uber Fouriershe Integrale,New York,1948.
%3.\textbf
\bibitem{EMOT}
{A. Erdelyi,  W. Magnus,  F. Oberhettingen,  and F. G.Tricomi},
\textit { Higher Transcendental functions.} New York, 1953.
%4.\textbf
\bibitem{GuRao}
{K.E.Gustatson and D.Rao}\textit{Numerical Range,}Springer,1996.
%5. \textbf
\bibitem{M.H}
{M.Haase,} \textit{The Functional Calculus for Sectorial Operators}
Operator Theory \textbf {169}, Birkhauser, 2006.
%6. \textbf
\bibitem{HaWi}{P.Hartman and A.Wintner}
\textit{On the Infinitesimal Generators of Integral Convolutions,}
Amer.UJ.Math.,64,273-298,1942.
%7. \textbf
\bibitem{HenTheo}{Hengartner W. and Theodorescu R.}\textit{Concentration
Functions}Academic Press, New York, 1973.
%8. \textbf
\bibitem{Ito}
{K. Ito}, \textit {On Stochastic Differential Equations.} Memoirs
Amer. Math. Soc. \textbf {4}, 1951.
%9. \textbf
\bibitem{Kac1}
{ M. Kac}, \textit {Probability and Related Topics in Physical Sciences}. Colorado, 1957.
%10.\textbf
\bibitem{Kac2}
{ M. Kac}, \textit {Some stochastic problems in physics and mathematics.} Dallas,1957.
%11. \textbf
\bibitem{KR}
{ M. G. Krein and  M. A. Rutman}, \textit {Linear Operators Leaving
Invariant a Cone
in a Banach Space.} Amer. Math. Soc., Translation \textbf {26}, 1950.
%12. \textbf
\bibitem{LEV}
{B.M.Levitan},\textit{Some Questions of the Theory of Almost Periodic Functions.}
Amer. Math. Soc. Transl.\textbf{28},1950.
%13. \textbf
\bibitem{RAD}
{J. Radon,} \textit{\"Uber lineare Functionaltransformationen and
Functionalgleichungen.} Sitzber. Akad. Wiss. Wien, 128, 1083-1121, 1919.
%13. \textbf
\bibitem{Sakh5}
{ L. A. Sakhnovich,} \textit {Integral Equations with Difference
Kernels.} Operator Theory \textbf {84}, Birkh\"auser, 1996.
%14. \textbf
\bibitem{Sakh6}
{L. A. Sakhnovich,} \textit {Levy Processes,I ntegral Equations, Statistical
 Physics: Connections and Interactions.}
 Operator Theory \textbf {225}, Birkh\"auser, 2012.
 
\bibitem{SakhArx}
L. A. Sakhnovich,
\textit{Convolution type form of the Ito representation of the infinitesimal generator for Levy processes,}
 arXiv:1212.3603
%15. \textbf
\bibitem{Sato}
{K.Sato},  \textit {Levy Processes and Infinitely Divisible
Distributions.} University
Press, Cambridge, 1999.
%16. \textbf
\bibitem{Sch}
{ W. Schoutens,} \textit {Levy Processes in Finance.} Wiley series
in
Probability and Statistics, 2003.
%17. \textbf
\bibitem{MSH}
{M.Sharpe,}\textit{Zeroes of Infinitely Divisible Densities,}
The Annals of Mathematical Statistics,Vol.40,No.4,1503-1505,1969.
%18. \textbf
\bibitem{STO}
{M.Stone,}\textit{Linear Transformation in Hilbert Space.} New York, 1932.
%19. \textbf
\bibitem{ThoBarn}
{M. Thomas and O. Barndorff} (ed.), \textit {Levy Processes:
Theory and Applications.} Birkh\"auser, 2001.
%20. \textbf
\bibitem{Titch}
{E. C. Titchmarsh}, \textit {Introduction to the Theory of Fourier Integrals.} Oxford, 1937.
%21. \textbf
\bibitem{Tuc}
{H.G.Tucker},\textit{The Supports of Infinitely Divisible Distribution
functions,} Proc. Amer. Math. Soc.,49, 436-440, 1975.
%22. \textbf
\bibitem{Zy}
{A.Zygmund},\textit{Trigonometric Series,} v.II, Cambridge, 1959.









\end{thebibliography}
\end{document}